\documentclass[12pt,reqno]{amsart}

\usepackage{amssymb}

\newcommand{\Z}{{\mathbb Z}}

\newcommand{\cE}{{\mathcal E}}
\newcommand{\cS}{{\mathcal S}}

\newcommand{\oA}{{\overline A}}
\newcommand{\oP}{{\overline P}}

\newcommand{\tA}{{\widetilde A}}
\newcommand{\tT}{{\widetilde T}}

\newcommand{\sig}{\sigma}

\newcommand{\Del}{\Delta}

\newcommand{\longc}{,\dotsc,}
\newcommand{\longp}{+\dotsb+}
\newcommand{\longe}{=\dotsb=}

\newcommand{\est}{\varnothing}

\newcommand{\seq}{\subseteq}

\newcommand{\stm}{\setminus}

\renewcommand{\(}{\left(}

\newcommand{\lfl}{\left\lfloor}
\newcommand{\rfl}{\right\rfloor}

\newcommand{\lcl}{\left\lceil}
\newcommand{\rcl}{\right\rceil}

\newcommand{\lfr}{\left\{}
\newcommand{\rfr}{\right\}}

\newcommand{\lpr}{\left(}
\newcommand{\rpr}{\right)}

\newtheorem{theorem}{Theorem}
\newtheorem{corollary}{Corollary}

\newcommand{\refc}[1]{\ref{c:#1}}
\newcommand{\reft}[1]{\ref{t:#1}}
\newcommand{\refs}[1]{\ref{s:#1}}
\newcommand{\refb}[1]{\cite{b:#1}}
\newcommand{\refe}[1]{\eqref{e:#1}}

\title{The structure of higher sumsets}
\author{Vsevolod F. Lev}
\email{seva@math.haifa.ac.il}
\address{Department of Mathematics, The University of Haifa at Oranim,
  Tivon 36006, Israel}
\makeatletter
\@namedef{subjclassname@2020}{\textup{2020} Mathematics Subject Classification}
\makeatother
\subjclass[2020]{Primary 11B13; Secondary 11D07, 11A99}
\keywords{Sumsets, Frobenius problem}

\begin{document}
\baselineskip=16pt

\begin{abstract}
Merging together a result of Nathanson from the early 70s and a recent
result of Granville and Walker, we show that for any finite set $A$ of
integers with $\min(A)=0$ and $\gcd(A)=1$ there exist two sets, the
``head'' and the ``tail'', such that if $m\ge\max(A)-|A|+2$, then the
$m$-fold sumset $mA$ consists of the union of the head, the appropriately
shifted tail, and a long block of consecutive integers separating them. We
give sharp estimates for the length of the block, and find all those sets
$A$ for which the bound $\max(A)-|A|+2$ cannot be substantially improved.
\end{abstract}

\maketitle

\section{Background, motivation, and summary of results.}

Let $A$ be a finite set of integers with $\min(A)=0$ and $\gcd(A)=1$. It is a
basic fact dating back to Frobenius and Sylvester that the additive semigroup
generated by $A$ contains all positive integers from some point on. The
largest integer that does not belong to this semigroup is called the
\emph{Frobenius number} of $A$. It is well-known that the Frobenius number of
the three-element set $A=\{0,a,b\}$ is $(a-1)(b-1)-1$, but no simple explicit
formula seems to exist for sets with four or more elements. The problem of
finding the Frobenius number is known as the \emph{linear diophantine problem
of Frobenius}; see~\refb{ra} for a comprehensive account.

The semigroup $\cS(A)$ generated by $A$ can be written as the infinite union
$\cS(A)=\{0\}\cup A\cup 2A\cup 3A\cup\ldots$\,, where $mA:=\{a_1\longp
a_m\colon a_1\longc a_m\in A\}$ are the \emph{sumsets} of $A$; that is, $mA$
is the set of all possible sums of $m$ elements of $A$, not necessarily
distinct. Notice, that as a result of $0\in A$, the sumsets satisfy
$\{0\}\seq A\seq 2A\seq 3A\seq\ldots$

In the light of the Frobenius-Sylvester observation, it is interesting to
investigate the structure of the individual sumsets $mA$. In this direction,
Nathanson has obtained the following nice result.
\begin{theorem}[Nathanson \refb{n}]\label{t:n}
Suppose that $A$ is a finite set of $n:=|A|\ge 3$ integers with $\min(A)=0$
and $\gcd(A)=1$, and let $l:=\max(A)$. Then there exist nonnegative integers
$h$ and $t$ and finite integer sets $H\seq[0,h-2]$ and $T\seq[0,t-2]$
depending only on $A$ such that for any integer
 $m\ge l^2(n-1)$ we have
\begin{equation}\label{e:nath}
  mA=H\cup[h,ml-t]\cup(ml-T).
\end{equation}
\end{theorem}
Loosely speaking, if $m$ is sufficiently large, then the sumset $mA$ consists
of a ``head'', an appropriately shifted ``tail'' (both head and tail
depending only on $A$ but not on $m$), and an interval separating them. The
fact that $h\le ml-t$, meaning that the interval is nonempty, is not stated
explicitly in~\refb{n}, but is implicit in the proof.

Both the Frobenius-Sylvester observation and the theorem of Nathanson
manifest the same phenomenon, which is that high-multiplicity sumsets of
properly normalized sets contain long blocks of consecutive integers.

Let $\cE(A)$ be the complement of the semigroup $\cS(A)$ in the set of all
positive integers. Following~\refb{gw}, we call the set $\cE(A)$ the
\emph{exceptional set} of $A$. (This set is also called the \emph{set of
gaps} of the semigroup $\cS(A)$.) The exceptional set is, therefore, the set
of all positive integers not representable as a nonnegative linear
combination of the elements of $A$, and its largest element $\max(\cE(A))$ is
the Frobenius number of $A$.

For brevity, we will write $\cS$ and $\cE$ instead of $\cS(A)$ and $\cE(A)$.

Keeping the notation of Theorem~\reft{n}, as an immediate consequence of the
theorem, the semigroup $\cS$ contains all integers larger than $h$, while an
integer $z\in[0,h]$ belongs to $\cS$ if and only if it belongs to $H$; that
is, $H=[0,h]\stm\cE$. Furthermore, from~\refe{nath} we get $m(l-A) = ml-mA =
T \cup [t,ml-h] \cup (ml-H)$. Therefore, denoting by $\cE'$ the exceptional
set of $l-A$, and repeating the argument above with $A$ replaced by $l-A$, we
obtain $T=[0,t]\stm\cE'$, and hence $ml-T=[ml-t,ml]\stm(ml-\cE')$. Thus,
Theorem~\reft{n} says that if $m\ge l^2(n-1)$, then $mA=[0,ml]\stm\lpr
\cE\cup(ml-\cE')\rpr$.

The bound $m\ge l^2(n-1)$ has been improved in a number of subsequent papers,
such as~\refb{wcc} or~\refb{gs}. Recently, Granville and Walker have shown
that $m\ge l-n+2$ suffices.
\begin{theorem}[Granville-Walker {\cite[Theorem~1]{b:gw}}]\label{t:gw}
Suppose that $A$ is a set of $n:=|A|\ge 3$ integers with $\min(A)=0$ and
$\gcd(A)=1$, and let $l:=\max(A)$. Then for any integer $m\ge M:=l-n+2$ we
have
\begin{equation}\label{e:gw}
  mA = [0,ml] \stm \big(\cE\cup(ml-\cE')\big),
\end{equation}
where $\cE$ and $\cE'$ are the exceptional sets of $A$ and $l-A$,
respectively.
\end{theorem}

As shown in~\refb{gs} and~\refb{gw}, the bound $m\ge l-n+2$ is best possible.
However, a natural question remains unanswered: can the sets $\cE$ and
$ml-\cE'$ overlap? Does the assumption $m\ge l-n+2$ guarantee the existence
of an interval separating these sets, and if so, what is the length of the
interval?

Our first goal here is to answer this question and give the Granville-Walker
result an alternative, surprisingly short proof. We show not only that
~\refe{gw} holds, but also that $\cE$ and $ml-\cE'$ are disjoint and,
moreover, there is a block of at least $(m-M+1)l+1$ consecutive integers
contained in $mA$ and separating $\cE$ from $ml-\cE'$.

\begin{theorem}\label{t:main}
Suppose that $A$ is a set of $n:=|A|\ge 3$ integers with $\min(A)=0$ and
$\gcd(A)=1$, and let $l:=\max(A)$. Then for any integer $m\ge M:=l-n+2$ we
have
  $$ mA = [0,ml] \stm \big(\cE\cup(ml-\cE')\big), $$
where $\cE$ and $\cE'$ are the exceptional sets of $A$ and $l-A$,
respectively. Moreover, the interval $[\max(\cE)+1,\min(ml-\cE')-1]$ is
nonempty and contained in $mA$; indeed, writing $l-1=k(n-2)+r$ with integers
$k\ge 1$ and $r\in[0,n-3]$, we have
  $$ \min(ml-\cE') - \max(\cE) \ge (m-M+1)l + \Del, $$
where
  $$ \Del = l(k-1)(n-3) + r k(n-3) + r^2 + k + 1 \ge 2. $$
\end{theorem}
We adopt the convention that if $\cE=\est$, then $\max(\cE)=-1$, and
similarly for $\max(\cE')$.

Notice that the quantity $\Del$ is ``normally'' much larger than $2$; for
instance, it is easy to show that $\Del>\(1-\frac1k-\frac1{n-2}\)l^2$.

The lower bound $(m-M+1)l+\Del$ established in the second part of the theorem
is best possible; equality is obtained, for instance, for the sets of the
form $A=\{0,d,2d\longc l\}\cup\{l-1\}$, where $d$ is a nontrivial divisor of
$l$. For these sets we have $n=\frac ld+2$ whence $l=(d-1)(n-2)+(n-3)+1$, so
that $k=d-1$ and $r=n-3$. Furthermore, the semigroup $\cS$ generated by $A$
is identical to that generated by the three-element set $\{0,d,l-1\}$, from
which, recalling the formula at the very beginning of this section, we derive
that $\max(\cE)=(d-1)(l-2)-1=k(l-2)-1$. Also, $\cE'=\est$ since $1\in(l-A)$.
Consequently, $\max(\cE)+\max(\cE')=k(l-2)-2$ and therefore
  $$ \min(ml-\cE') - \max(\cE) = ml - k(l-2)+2 = (m-M+1)l + \Del, $$
where the second equality can be verified by cancelling out the common
summand $ml$ and expressing the rest in terms of the parameters $n$ and $d$.

More generally, one can consider sets of the form
 $$ A:=\{0,d,2d\longc sd\}\cup\{sd-1-td,sd-1-(t-1)d\longc sd-1\}, $$
where $s,d$ and $t$ are positive integers satisfying $t<s$ and $2\le
d<\frac{s}{t}+1$. We have $l=sd$, $n=s+t+2$, and the exceptional set of $A$
is the same as that of the set $\{0,d,sd-1-td\}$. It follows that
$\max(\cE)=(d-1)((s-t)d-2)-1$, and consequently,
  $$ \min(ml-\cE') - \max(\cE) = ml - (d-1)((s-t)d-2)+2 = (m-M+1)l + \Del. $$
(The last equality can be verified by observing that $k=d-1$, as it follows
from the assumption $d<\frac st+1$.)

We prove Theorem~\reft{main} in Section~\refs{main}. The proof is a further
elaboration on the ideas from \cite{b:gs,b:gw}.

\smallskip
Our second goal is to investigate the corresponding stability problem and
determine those sets $A$ such that for~\refe{gw} to hold, one needs $m$ to be
almost as large as $l-n+2$. Some results in this direction are obtained
in~\refb{gw} where the sets requiring $m\ge l-n$ are fully described (under
the technical assumptions $l\ge 9$ and $n\ge 5$). Citing from~\refb{gw},

\smallskip
\begin{quotation}
\emph{Indeed, our proofs are sufficiently flexible that one can go on and
prove that $mA=[0,ml]\stm\big(\cE\cup(ml-\cE')\big)$ holds for all
$m\ge\max\{1,l-n-C\}$ for even larger values of $C$, except in some explicit
finite set of families of sets $A$, though the number of cases seems to grow
prohibitively with $C$.}
\end{quotation}

\smallskip
We prove the following result showing that the ``tough'' sets are,
essentially, dense subsets of the set $\{0,1\}\cup[m+2,l]$ or of its mirror
reflection $[0,l-(m+2)]\cup\{l-1,l\}$.

\begin{theorem}\label{t:stab}
Let $A$ be a set of $n:=|A|\ge 6$ integers with $\min(A)=0$ and $\gcd(A)=1$.
Write $l:=\max(A)$, and let $\cE$ and $\cE'$ be the exceptional sets of $A$
and $l-A$, respectively. Then for any integer $m$ with
\begin{equation}\label{e:Mstab}
  m \ge \max \lfr l-\frac32\,n+\frac92,\,\frac23\,(l-n+2) \rfr
\end{equation}
we have
\begin{equation}\label{e:stab1}
   mA = [0,ml] \stm \big(\cE\cup(ml-\cE')\big),
\end{equation}
except if either $\{0,1\}\seq A\seq\{0,1\}\cup[m+2,l]$, or
 $\{l-1,l\}\seq A\seq[0,l-(m+2)]\cup\{l-1,l\}$, in which cases~\refe{stab1}
does not hold.
\end{theorem}

Notice that each of the sets $\{0,1\}\cup[m+2,l]$ and
$[0,l-(m+2)]\cup\{l-1,l\}$ has size $l-m+1<\frac32\,n$, the inequality
following from~\refe{Mstab}; this shows that their $n$-element subsets are
contained therein with density exceeding $2/3$.

We also remark that for $n\ge 6$, the maximum in the right-hand side
of~\refe{Mstab} is always smaller than $l-n+2$, the bound of
Theorem~\reft{main}.

The proof of Theorem~\reft{stab} is presented in Section~\refs{stab}.

\section{The toolbox}\label{s:toolbox}

In this section we collect various results used in the proofs of
Theorems~\reft{main} and~\reft{stab}.

We use the standard set addition notation: if $B$ and $C$ are subsets of an
additively written group, then the sumset $B+C$ is defined to be the set
$\{b+c\colon b\in B,\,c\in C\}$.

\begin{theorem}[Olson {\cite[Theorem~1]{b:o}}]\label{t:olson}
Suppose that $A$ and $B$ are finite, nonempty subsets of a group. If $0\in
A$, then either $A$ is contained in the subgroup of all those group elements
$z$ satisfying $(A+B)+z=A+B$, or $|A+B|\ge|B|+\frac12|A|$.
\end{theorem}

Note that if $A$ is not contained in a proper subgroup, then the first
alternative is ruled out unless $A+B$ is the whole group. With this
observation in mind, arguing inductively we obtain the following corollary.
\begin{corollary}\label{c:olson32}
Suppose that $A$ is a finite subset of an abelian group $G$, and $m\ge 1$ is
an integer. If $A$ is not contained in a proper coset, then either $mA=G$, or
$|mA|\ge\frac{m+1}2\,|A|$.
\end{corollary}

Here is yet another immediate corollary of Olson's theorem.
\begin{corollary}\label{c:BF}
If $A$ and $B$ are finite, nonempty subsets of an abelian group $G$ with
$|A+B|\le|B|+1$, $A+B\ne G$, and $|A|\ge 3$, then $A$ is contained in a coset
of a proper subgroup of $G$.
\end{corollary}

\begin{theorem}[Scherk \refb{s}]\label{t:scherk}
Suppose that $B$ and $C$ are finite subsets of an abelian group with $0\in
B\cap C$. If $0$ has a unique representation in $B+C$ (which is $0=0+0$),
then $|B+C|\ge|B|+|C|-1$.
\end{theorem}

Iterating Theorem~\reft{scherk}, we get
\begin{corollary}[Alon {\cite[Corollary~2.3]{b:a}}]\label{c:a}
Suppose that $m\ge 1$ is an integer, and $A$ is a finite subset of an abelian
group. If $0\notin mA$, then $|A\cup(2A)\cup\dotsb\cup(mA))|\ge m|A|$.
\end{corollary}

Next, we need a result often referred to as \emph{Freiman's
$(3n-3)$-theorem}.
\begin{theorem}[Freiman \cite{b:f}]\label{t:3n-3}
Suppose that $A$ is a set of $n:=|A|\ge 3$ integers with $\min(A)=0$ and
$\gcd(A)=1$, and let $l:=\max(A)$. Then
  $$ |2A| \ge \min\{l,2n-3\} + n. $$
\end{theorem}

A sequence of elements of an abelian group is called \emph{zero-sum-free} if
it does not have a finite, nonempty subsequence with the zero sum of its
terms. The following theorem describing the structure of long zero-sum-free
sequences in finite cyclic groups was proved by Savchev and Chen and,
simultaneously and independently, by Yuan.
\begin{theorem}%
[{Savchev-Chen \cite[Theorem~8]{b:sc}, Yuan \cite[Theorem~3.1]{b:y}}]%
  \label{t:sc}
Suppose that $u$ and $l$ are positive integers, and that $(a_1\longc a_u)$ is
a zero-sum-free sequence of elements of the cyclic group of order $l$. If
$u>l/2$, then there exist positive integers $x_1\longc x_u$ with $x_1\longp
x_u<l$ and a group element $a$ of order $l$ such that $a_i=x_ia$ for any
$i\in[1,n]$.
\end{theorem}

Finally, we list two theorems by the present author.

\begin{theorem}[{\cite[Proposition~1]{b:ast}}]\label{t:ast}
Suppose that $X=(x_1\longc x_u)$ is a nonempty sequence of positive integers
written in an increasing order: $1\le x_1\le\dotsb\le x_u$. If $X$ has fewer
than $2u$ distinct subsequence sums, then $x_2\longc x_u$ are all divisible
by $x_1$, and $x_{i+1}\le x_1\longp x_i$ for each $i\in[1,u-1]$.
\end{theorem}

\begin{theorem}[{\cite[Theorems~1 and 3 (ii)]{b:lopt}}]\label{t:opt}
Suppose that $A$ is a set of $n:=|A|\ge 3$ integers with $\min(A)=0$ and
$\gcd(A)=1$. Let $l:=\max(A)$ and write $l-1=k(n-2)+r$ with $k\ge 1$ and
$r\in[0,n-3]$ integers.
\begin{itemize}
\item[(i)] If $m\ge 2k$, then $[kl-k(n-1-r),(m-k)l+k(n-1-r)]\seq mA$;
\item[(ii)] if $m\ge 3k$, then indeed $mA$ contains a block of at least
    $(m-k)l+k(n-1-r)+1$ consecutive integers.
\end{itemize}
\end{theorem}

An immediate corollary of Theorem~\reft{opt} is an estimate for the Frobenius
number of a set of given size and ``diameter'', originally proved by Dixmier.
\begin{corollary}[Dixmier {\cite[Theorem~3]{b:d}}]\label{c:Dix}
Suppose that $A$ is a set of $n:=|A|\ge 3$ integers with $\min(A)=0$ and
$\gcd(A)=1$. Let $l:=\max(A)$ and write $l-1=k(n-2)+r$ with integer $k$ and
$0\le r<n-2$. Then $\max(\cE)\le k(l-n+r+1)-1$, where $\cE$ is the
exceptional set of $A$.
\end{corollary}

\section{Proof of Theorem~\reft{main}}\label{s:main}

Recall that we have defined integers $k\ge 1$ and $r\in[0,n-3]$ by
$l-1=k(n-2)+r$.

Our first goal is to show that if $n\ge 4$, then $[0,kl]\stm\cE\seq mA$. (We
remark that the assumption $n\ge 4$ is essential and cannot be dropped; say,
if $A=\{0,1,l\}$, then $M=k=l-1$, $\cE=\est$, and $[0,kl]\nsubseteq MA$ as,
for instance, $l(l-1)-1\notin MA$.)

Suppose, for a contradiction, that $g\in[0,kl]\stm\cE$ is an integer not
representable as a sum of $m$ or fewer elements of $A$. If $g$ is a multiple
of $l$, then we can write $g=vl$ with an integer $v$ satisfying $v>m$ in view
of $g\notin mA$ and $l\in A$. On the other hand, $v\le k$ since $g\le kl$.
Hence,
  $$ k \ge v > m \ge M = l-(n-2)
                            \ge (k(n-2)+1) - (n-2) = (k-1)(n-2) +1 \ge k, $$
a contradiction. Thus, $g$ is not a multiple of $l$.

Since $g\notin\cE$, we can write
\begin{equation}\label{e:repg}
  g = a_1\longp a_u + lv
\end{equation}
where $u,v\ge 0$ and $a_1\longc a_u\in A\stm\{0,l\}$ are integers; indeed,
since $g$ is shown above not to be divisible by $l$, we have $u\ge 1$ and
$v\le k-1$. Suppose that the representation~\refe{repg} has the smallest
value of the parameter $u$ possible, among all representations of $g$ in this
form. We show that in this case $u+v\le l-(v+1)(n-2)$; since the quantity in
the right-hand side is at most $l-n+2=M$, this will imply $g\in MA\seq mA$,
contradicting the choice of $g$.

We notice that, in view of $u+v>m\ge M=l-n+2$ (which follows from the
assumption $g\notin mA$) and $v\le k-1$,
\begin{multline}\label{e:uv2}
  u \ge l-n+3 - v \ge (k(n-2)+r+1) - (n-2) + 1 -(k-1) \\
            \ge (k-1)(n-3)+2 \ge k+1 \ge v+2.
\end{multline}

An important observation originating from~\cite{b:gs,b:gw} is that the
sequence $(a_1\longc a_u)$ reduced modulo $l$ is zero-sum-free; that is, does
not have any nonempty subsequences with the sum of their elements divisible
by $l$. Indeed, if we had, say,
 $a_1\longp a_s=wl$ with $s\in[1,u]$ and $w$ integers, this would lead to
$g=a_{s+1}\longp a_u+(v+w)l$, contradicting minimality of $u$.

As a result of $(a_1\longc a_u)$ being zero-sum-free, the $u-v-1$ sums
  $$ \sig_s:=a_1\longp a_s,\quad v+2\le s\le u $$
are pairwise distinct modulo $l$, cf.~\refe{uv2}. Following \refb{gw}, we
claim that, moreover, these sums are also distinct modulo $l$ from all
elements of the sumset $(v+1)A$. To see this, suppose that $\sig_s=b_1\longp
b_t+wl$ with integers $v+2\le s\le u$, $1\le t\le v+1$, and $w$, and elements
$ b_1\longc b_t\in A\stm\{0,l\}$. Then $\sig_s>0$ implies $w\ge -(t-1)\ge
-v$, showing that $v+w\ge 0$. Consequently, $g= b_1\longp b_t+a_{s+1}\longp
a_u+(v+w)l$ is a representation of $g$ contradicting minimality of $u$ in
view of $t\le v+1<v+2\le s$.

Let $\oA$ denote the canonical image of $A$ in the quotient group $\Z/l\Z$.
As we have just shown, the sumset $(v+1)\oA$ is disjoint from the set of all
$u-v-1$ sums $\sig_s$, $s\in[v+2,u]$, taken modulo $l$. It follows that
\begin{equation}\label{e:vA}
  |(v+1)\oA|\le l-(u-v-1).
\end{equation}
Since $u+v\ge m+1\ge M+1=l-n+3$, this gives
\begin{multline*}
  |(v+1)\oA| \le l+v+1-u
             \le (l + v + 1) - (l-n+3-v)
             = 2v+n-2 = |\oA|+2v-1.
\end{multline*}
As a result, $v\ge 1$, and there exists an integer $v_0\in[1,v]$ with
$|(v_0+1)\oA|\le|v_0\oA|+1$. By Corollary~\refc{BF}, the set $\oA$ is
contained in a coset of a proper subgroup. Therefore, the original set $A$ is
contained in an arithmetic progression with the endpoints $0$ and $l$ and the
difference larger than $1$, contradicting the assumption $\gcd(A)=1$.

We have thus shown that if $n\ge 4$, then
  $$ [0,kl]\stm\cE\seq mA. $$
Switching the roles of $A$ and $l-A$, we conclude that $[0,kl]\stm\cE'\seq
m(l-A)$; equivalently,
  $$ [(m-k)l,ml]\stm(ml-\cE')\seq mA. $$
If $m\le 2k$, then the intervals $[0,kl]$ and $[(m-k)l,ml]$ jointly cover the
whole interval $[0,ml]$, and we immediately obtain~\refe{gw}. If $m>2k$, then
by Theorem~\reft{opt}~(i) we have $[kl,(m-k)l]\seq mA$ which, again, leads
to~\refe{gw}.

This proves the first assertion of the theorem in the case where $n\ge 4$.
Since the case $n=3$ has received a nice and simple dedicated treatment
in~\cite[Theorem~4]{b:gs}, we ignore it here and proceed to the second
assertion (without assuming $n\ge 4$ any longer).

Write $e:=\max(\cE)$ and $e':=\max(\cE')$. By Theorem~\reft{opt}~(ii), the
sumset $(3k)A$ contains a block, say $B$, of $|B|=2k l+k(n-1-r)+1$
consecutive integers. Since $|B|>l$, while the difference between any two
consecutive elements of $\cE$ is easily seen not to exceed $l$ (and indeed,
not to exceed the smallest nonzero element of $A$) and similarly for $\cE'$,
we have $B\seq[e+1,3kl-e'-1]$. Therefore $3kl-e'-e-1\ge|B|=2kl+k(n-1-r)+1$
whence $kl-e-e'-1\ge k(n-1-r)+1$. It follows that
\begin{align*}
  \min(ml-\cE')-\max(\cE)
     &= ml-e-e' \\
     &\ge ml -kl + k(n-1-r)+2 \\
     &= (m-M+1)l + \Del
\end{align*}
(verifying the last equality is tedious, but straightforward, and we omit the
computation).

\section{Proof of Theorem~\reft{stab}}\label{s:stab}
We keep writing $l-1=k(n-2)+r$ with integers $k\ge 1$ and $r\in[0,n-3]$.
Notice that if $k\ge 2$, then from~\refe{Mstab} and the assumption $n\ge 6$,
we have
  $$ m \ge l-\frac32n+\frac92 \ge k(n-2) - \frac32\,(n-2) + \frac52
       = \lpr k-\frac32\rpr (n-2) + \frac52 \ge 4k-\frac72 $$
whence
\begin{equation}\label{e:m4}
  m\ge 4k-3 \ge k.
\end{equation}
Clearly, this resulting estimate remains valid also if $k=1$.

We split the proof into several parts.

\subsection{Sufficiency}
First, we show that $\{0,1\}\seq A\seq\{0,1\}\cup[m+2,l]$ implies
 $mA\ne[0,ml]\stm(\cE\cup(ml-\cE'))$; by symmetry, the same conclusion
follows also from $\{l-1,l\}\seq A\seq[0,l-(m+2)]\cup\{l-1,l\}$. We actually
show that $m+1\notin\cE$ and $m+1\notin(ml-\cE')$; this proves the assertion
since, clearly, $m+1\notin mA$ while $m+1\in[0,ml]$.

The relation $m+1\notin\cE$ is, indeed, trivial since $\cE=\est$ in view of
$1\in A$. To prove that $m+1\notin(ml-\cE')$ we rewrite this as
$ml-m-1\notin\cE'$ and notice that by Corollary~\refc{Dix}, we have
 $\max(\cE')\le k(l-n+r+1)-1$. Thus, it
suffices to show that $ml-m-1>k(l-n+r+1)-1$; that is, $m(l-1)>k(l-n+r+1)$.
However, the last inequality is clearly true in view of~\refe{m4} and since
 $r\le n-3$.

The rest of the argument deals with necessity; we want to show that, under
the stated assumptions, we have $mA=[0,ml]\stm\big(\cE\cup(ml-\cE')\big)$,
unless either $\{0,1\}\seq A\seq\{0,1\}\cup[m+2,l]$, or
 $\{l-1,l\}\seq A\seq[0,l-(m+2)]\cup\{l-1,l\}$.

\subsection{The general setup}
The case where $l\le n$ is easy to analyze, and we assume below that $l>n$.

We show that if $m$ satisfies~\refe{Mstab}, then either
 $1\in A\seq\{0,1\}\cup[m+2,l]$, or
\begin{equation*}\label{e:3008a}
  [0,kl] \stm \cE\seq mA.
\end{equation*}
Once this is established, the proof can be easily completed. Namely, applying
the assertion to the set $l-A$, we conclude that either
 $(l-1)\in A\seq[0,l-(m+2)]\cup\{l-1,l\}$, or
\begin{equation*}\label{e:3008b}
  [(m-k)l,ml] \stm (ml-\cE') \seq mA.
\end{equation*}
Taking into account that, by Theorem~\reft{opt} (i), if $m\ge 2k$, then
$[kl,(m-k)l]\seq mA$, we conclude that
  $$ [0,ml] \stm \big(\cE\cup(ml-\cE')\big) \seq mA, $$
and the converse inclusion is trivial.

We thus assume that $g\in[0,kl]\stm\cE$ is an integer with $g\notin mA$, and
show that $A\cap[2,m+1]=\est$ and that $1\in A$.

\subsection{Preliminaries}

From~\refe{m4} we have $kA\seq mA$. Therefore $g$ is not a multiple of $l$,
as if we had $l\mid g$, then $g/l\le k$ and $l\in A$ would lead to
 $g=(g/l)l\in kA\seq mA$.

Since $g\notin\cE$, we can write
\begin{equation}\label{e:repgen}
  g=a_1\longp a_u+vl
\end{equation}
where $u,v\ge 0$ and $a_1\longc a_u\in A\stm\{0,l\}$ are integers. Speaking
about \emph{representations} of $g$ we always mean representations of this
particular form. In view of $g\notin mA$, for any representation we have
\begin{equation}\label{e:uvlarge}
  u+v\ge m+1.
\end{equation}
Since $g$ is not a multiple of $l$, we also have $u>0$ and $v\le k-1$. Hence,
recalling~\refe{m4},
  $$ u + k-1 \ge u+v \ge m+1 \ge 4k-2; $$
therefore, $u\ge 3k-1\ge v+2k$ and
\begin{equation}\label{e:3008c}
  u+v\ge 4k-2.
\end{equation}

Without loss of generality, we assume that the representation~\refe{repgen}
has the smallest possible value of the parameter $u$ among all
representations of $g$; this assumption will be referred to as
\emph{minimality} (of $u$).

\subsection{The zero-sum-free property}
Repeating literally the argument from the proof of Theorem~\reft{main}, we
conclude that the sequence $(a_1\longc a_u)$ is zero-sum-free modulo $l$, the
sums
  $$ \sig_s:=a_1\longp a_s,\quad v+2\le s\le u, $$
taken modulo $l$ are distinct from each other and from all elements of the
sumset $(v+1)\oA$, where $\oA:=A\pmod l$, and as a result,~\refe{vA} holds
true.

On the other hand, by Corollary~\refc{olson32}, since $\oA$ is not contained
in a proper subgroup (as it follows from $\gcd(A)=1$), we have
$|(v+1)\oA|\ge\frac{v+2}2\,|\oA|$. Therefore,
  $$ l-u+v+1 \ge \frac{v+2}2\,(n-1); $$
consequently, since $u+v\ge m+1\ge l-\frac32\,n+\frac{11}2$ by~\refe{uvlarge}
and~\refe{Mstab},
\begin{align*}
  l &\ge u-v-1+\frac{v+2}2\,(n-1) \\
    &= (u+v) + \frac v2\,(n-5) + n-2 \\
    &\ge l-\frac32\,n + \frac{11}2 + \frac v2\,(n-5) + n-2 \\
    &= l + \frac12(v-1)(n-5) + 1
\end{align*}
which is wrong whenever $v\ge 1$ and $n\ge 6$. Thus, $v=0$. As a result,
by~\refe{uvlarge} and~\refe{Mstab},
\begin{equation}\label{e:ul1}
   u \ge l-\frac32\,n+\frac{11}2 \ \text{and}\ u \ge \frac23\,(l-n+2)+1,
\end{equation}
and by~\refe{3008c}
\begin{equation}\label{e:ul2}
  u \ge 4k-2.
\end{equation}

On the other hand, we have $u<l$ due to the basic fact that a zero-sum-free
sequence in a cyclic group has length smaller than the size of the group.

We investigate separately three possible cases: $g<l$, $l<g<2l$, and $g>2l$.

\subsection{The case where $g<l$}\label{s:g<l}
Consider the sets
  $$ D_j := (jA\stm(j-1)A)\cap[0,l],\quad j\in[1,u], $$
where $0A=\{0\}$ is assumed; thus, $D_1=A\stm\{0\}$, $g\in D_u$, and $D_j\seq
jA$. The last relation shows that every element $d\in D_j$ can be written as
$d=a+b$ with $a\in A$ and $b\in (j-1)A$; moreover, if $j\ge 2$, then
$b\notin(j-2)A$ as otherwise we would have $d=a+b\in(j-1)A$. Consequently,
$b\in D_{j-1}$ and we conclude that $D_j\seq D_{j-1}+A$. This shows, in
particular, that if $D_{j-1}=\est$, then also $D_j=\est$; therefore, from
$g\in D_u$, all sets $D_1\longc D_u$ are nonempty. Since these sets are
pairwise disjoint and contained in the interval $[0,l]$, using~\refe{ul1} and
the assumption $n<l$ (made at the very beginning of the proof) we obtain
$u\ge 3$, and then, using~\refe{ul1} once again,
  $$ |D_2|\longp|D_{u-1}|
     \le l+1 - (1+|D_1| + |D_u|)
     \le l - n
     <   2(u-2). $$
As a result, there is an integer $j\in[2,u-1]$ with $|D_j|=1$.

We now observe that for any subset $I\seq[1,u]$ of size $h:=|I|$ we have
$\sum_{i\in I}a_i\in D_h$ as if the sum were lying in $hA\stm D_h$, then it
would be representable as a sum of fewer than $h$ elements of $A$, leading to
a representation of $g$ with fewer than $u$ summands.

It follows that all sums $\sum_{i\in I}a_i$ with $|I|=j$ are equal to the
same number, which is the element of $D_j$. As a result, all elements
$a_1\longc a_u$ are equal to each other, and we denote by $a$ their common
value; thus, $g=ua\notin(u-1)A$ by the minimality of $u$.

By the box principle, from $g\notin 2A$ it follows that $|A\cap[0,g]|\le \lcl
g/2\rcl$ whence
  $$ n = |A| \le |A\cap[0,g]| + (l-g) \le l - \lfl g/2\rfl. $$
Hence, $\lfl g/2\rfl\le l-n$ showing that $au=g\le 2(l-n)+1 < 3u$, as it
follows from~\refe{ul1}. Therefore, $a<3$; that is, $a=1$ or $a=2$.

If $a=1$, then $g=u\ge m+1$ and $A\cap[2,u]=\est$ in view of $1=a\in A$ and
by the minimality of $u$, completing the proof. If $a=2$, then $g=2u$ and
$2=a\in A$ along with the minimality of $u$ is easily seen to imply
$A\cap[3,u]=\est$. Moreover, $A$ is disjoint from the set of all even
integers in the range $[u+1,2u]$. A simple counting now gives
  $$ |A|\le 3+\lfl\frac u2\rfl+(l-2u) \le l-\frac32\,u+3<n $$
(following from~\refe{ul1}), a contradiction.

\subsection{The case where $l<g<2l$}
In this case in view of $g<kl$ we have $k\ge 2$; consequently, $u\ge 4k-2\ge
6$ by~\refe{ul2}.

Modifying slightly the argument employed in the case $g<l$, this time we
define
  $$ D_j := (jA\stm(j-1)A)\cap[0,2l],\quad j\in[1,u]. $$
Since $k\ge 2$, we have $l\ge 2n-3$; hence, $|2A|\ge 3n-3$ by
Theorem~\reft{3n-3}. Therefore, arguing as in the case $g<l$,
\begin{align*}
  |D_3|\longp|D_{u-1}|
     &\le 2l+1 - (1+|D_1|+|D_2|) - |D_{u}| \\
     &\le   2l + 1 - |2A| - 1 \\
     &\le 2l     - (3n-3)  \\
     &< 2(u-3),
\end{align*}
the last inequality following from~\refe{ul1}. Consequently, there is an
integer $j\in[3,u-1]$ with $|D_j|=1$. As above, this leads to $a_1\longe
a_u$, and denoting this common value by $a$, we have
$g=ua\notin(u-1)A$. 

Since $u\ge 6$, we have $g\notin4A$. By the box principle,
$|2A\cap[0,g]|\le\lcl g/2\rcl$ implying
 $|2A|\le\lcl g/2\rcl + (2l - g) = 2l - \lfl g/2\rfl$ and then, by
Theorem~\reft{3n-3} and~\refe{ul1},
\begin{gather*}
  \lfl g/2 \rfl \le 2l-3n+3, \\
  g-1 \le 2(2l-3n+3), \\
  g\le 4\lpr l-\frac32\,n+\frac74 \rpr < 4u.
\end{gather*}
Since $g=ua$, we conclude that $a\le 3$. On the other hand, $a\ge 2$ in view
of $al>ua=g>l$.

If $a=2$, then $g=2u$, $2=a\in A$, and $l/2<u<l$. Consequently, from the
minimality of $u$ we derive that $A\cap[3,u]=\est$, and $A$ does not contain
any even integers in the range $[u+1,l]$. Therefore
  $$ |A| \le 3 + \lfl \frac{l-1}{2} \rfl - \lfl \frac u2 \rfl
       \le 3 + \frac l2 - \frac u2 = \frac12\,(l-u+6) < n, $$
a contradiction.

Assume now that $a=3$. If $A$ contained an element in the interval $[4,u]$,
then denoting this element by $b$, we could replace $b$ summands in the
representation $g=3\longp 3$ with just three summands $b+b+b$, contradicting
minimality of $u$; therefore $A\cap[4,u]=\est$.

Let
  $$ B_i:=\{b\in A\colon u+1\le b\le l
                     \ \text{and}\ b\equiv i\!\!\!\pmod 3\},\ i\in[0,2]. $$
By minimality of $u$ we have $B_0=\est$; in particular,
 $l\not\equiv 0\pmod 3$. Furthermore, $|A\cap[0,3]|\le 3$: otherwise
considering $j\in\{1,2\}\seq A$ such that $l+j\equiv 0\pmod 3$ we get a
contradiction. Next, if one of $B_1$ and $B_2$ is empty, then by~\refe{ul1}
  $$ n = |A| \le 3 + \frac13\,(l-u+2)
                         \le 3 + \frac13\,\lpr\frac32\,n-\frac72 \rpr < n, $$
a contradiction. Finally, if both $B_1$ and $B_2$ are nonempty, then the
sumset $B_1+B_2$ consists of multiples of $3$, and hence its smallest element
exceeds $g=3u$. Since $\min(B_i)\le l-3(|B_i|-1),\ i\in\{1,2\}$, with at
least one of these two inequalities strict, and since $B_1,B_2\ne\est$
implies $A\cap\{1,2\}=\est$, we have in this case
\begin{multline*}
  3u < (l-3(|B_1|-1)) + (l - 3(|B_2|-1)) - 1 \\
                                      = 2l-3(|B_1|+|B_2|)+5 = 2l-3(n-2)+5
\end{multline*}
whence $3u\le 2(l-n)+5-(n-6)\le 2(l-n)+5$, contradicting~\refe{ul1}.

\subsection{The case where $g>2l$}
Since $g<kl$, in this case we have $k\ge 3$.

As shown above, the sequence $(a_1\longc a_u)$ reduced modulo $l$ is
zero-sum-free. From~\refe{ul1}, and since $k\ge 3$, the length of this
sequence is
\begin{equation}\label{e:ularge}
   u \ge l - \frac32\,(n-2)+\frac52
                 \ge l - \frac32\,\frac{l-1}k + \frac52 \ge \frac l2+3.
\end{equation}
Consequently, by Theorem~\reft{sc}, there are integers $a\in[1,l-1]$ and
$x_1\longc x_u\ge 1$ such that $\gcd(a,l)=1$, $x_1\longp x_u<l$, and
$a_i\equiv ax_i\pmod l$ for each $i\in[1,u]$. As an immediate corollary,
$u<l$. Renumbering, we can assume that $x_1\le\dotsb\le x_u$.

For $j\in\{u-1,u\}$, let $P_j$ denote the set of all nonempty subsequence
sums of the sequence $(a_i)_{i\in[1,u]\stm\{j\}}$. The elements of the set
$a_j+P_j$ are distinct modulo $l$ from the elements of $A$: if, say,
$a_j+\sum_{i\in I}a_i=b+wl$ with a nonempty set $I\seq[1,u]\stm\{j\}$, and
with integers $b\in A\stm\{0,l\}$ and $w$, then $b+wl>0$ whence $w\ge 0$ and
$g=b+\sum_{i\in[1,u]\stm(I\cup\{j\})} a_i+wl$, contradicting the minimality
of $u$. Recalling the notation $\oA=A\pmod l$, and writing $\oP_j:=P_j\pmod
l$, we conclude that
  $$ |\oP_j|\le l-|\oA| = l-n+1\le 2(u-1)-2, $$
the last inequality being a consequence of~\refe{ul1}.
Since $x_1\longp x_u<l$, reducing the set $P_j$ modulo $l$ does not affect
its size.
Therefore, the subsequence sum set of the integer sequence
$(x_i)_{i\in[1,u]\stm\{j\}}$ has size $|P_j|+1=|\oP_j|+1\le 2(u-1)-1$.
Applying Theorem~\reft{ast} twice, first time with $j=u-1$, and the second
time with $j=u$, we conclude that all integers $x_1\longc x_u$ are divisible
by $x_1$, and that $x_{i+1}\le x_1\longp x_i$ for each $i\in[1,u-1]$.

Assuming for a moment that $x_1\longc x_u$ are not all equal to each other,
let $s\in[1,u-1]$ be the smallest integer such that $x_{s+1}>x_s$. Then the
subsequence sum set of the truncated sequence $(x_1\longc x_s)$ is the
arithmetic progression $\{0,x_1,2x_1\longc sx_1\}$; therefore, there exists
an index set $I\seq[1,s]$ with $|I|\ge 2$ such that $x_{s+1}=\sum_{i\in I}
x_i$. It follows that
  $$ \sum_{i\in I} a_i = a_{s+1} + wl $$
where $w$ is an integer and $w\ge 0$ in view of $wl+a_{s+1}>0$. This,
however, is easily seen to contradict the minimality of $u$.

We have thus shown that $x_1\longc x_u$ are all equal to each other; hence,
$a_1\longc a_u$ are equal to each other, too, and we denote their common
value by $a$. We notice that
\begin{equation}\label{e:2k-1}
  3 \le a\le 2k-1:
\end{equation}
the lower bound follows from $2l<g=ua<al$, and the upper bound from
$kl>g=ua>\frac 12\,al$, cf.~\refe{ularge}.

Comparing~\refe{2k-1} and~\refe{ul2}, we conclude that $u\ge 2a$.
Furthermore, we have $A\cap[a+1,u]=\est$, as if there existed an integer
$b\in A\cap[a+1,u]$, then we could reduce the number of summands in the
representation~\refe{repgen} by replacing the $b$-term sum $a\longp a$ with
the $a$-term sum $b\longp b$. We thus can partition $A$ as $A=A_0\cup A_1$
where $A_0:=A\cap[0,a]$ and $A_1=A\cap[u+1,l]$.

Let $H:=\lcl g/l\rcl-1$, and suppose that $1\le h\le H$. If
 $za\in A_0+hA_1$ with an integer $z$, then $z\ge h+2$ in view of
  $$ (h+1)a \le \frac12\,(h+1)u \le uh < \min(A_0+hA_1), $$
and similarly, $z\le u$ in view of
  $$ (u+1)a = a+g > a+Hl \ge a+hl = \max(A_0+hA_1). $$
Therefore, $h+2\le z\le u$, which is easily seen to contradict the minimality
of $u$ (replace $z$ summands of the representation $g=a\longp a$ with the
$(h+1)$-summand representation of $za$). Consequently, the sumset $A_0+hA_1$
does not contain multiples of $a$. Since this holds for any $h=1\longc H$,
letting
 $T:=A_1\cup\dotsb\cup HA_1$, we conclude that, indeed, $A_0+T$ does not
contain multiples of $a$; that is, $0\notin\tA_0+\tT$ where we let
$\tA_i:=A_i\pmod a\ (i=0,1)$ and $\tT:=T\pmod a$.

By the box principle, from $0\notin\tA_0+\tT$ it follows that
 $|\tT|\le a-|\tA_0|=a-|A_0|+1$. The fact that $0\notin\tA_0+\tT$ also
implies $0\notin\tT$; consequently, $|\tT|\ge H|\tA_1|$ by
Corollary~\refc{a}. Comparing these estimates, we obtain
\begin{equation}\label{e:0509a}
  H|\tA_1|\le a-|A_0|+1.
\end{equation}
On the other hand, by~\refe{ularge},
\begin{equation}\label{e:0509b}
  H > \frac gl\, - 1 = \frac ul\,a - 1 > \frac12\,a - 1
\end{equation}
whence, using~\refe{2k-1},
\begin{equation}\label{e:0509c}
  3H \ge \frac32\,a-\frac32 > a-|A_0|+1.
\end{equation}
From~\refe{0509a} and~\refe{0509c} we obtain $|\tA_1|\le 2$.

If $|\tA_1|=1$, then $A_1$ is contained in one single residue class modulo
$a$. Hence, by~\refe{ul1} and in view of~\refe{2k-1},
  $$ |A_1| \le \frac{l-u-1}a + 1
         \le \frac1a\,\Big( \frac32\,n-\frac{13}2\Big) + 1
                                          \le \frac12\,n-\frac{7}6. $$
On the other hand, from~\refe{0509a}
and~\refe{ularge},
  $$ a-|A_0|+1 \ge H > \frac gl-1 = \frac ul\,a-1
                                       > \Big( 1-\frac3{2k}\Big) \,a-1 $$
whence, by~\refe{2k-1},
  $$ |A_0| < \frac3{2k}\,a+2 < 5. $$

Combining the upper bounds, we get
  $$ n = |A_0|+|A_1| \le 4 + \frac12\,n-\frac{7}6 $$
which is wrong for $n\ge 6$.

Thus $|\tA_1|=2$. Substituting back to~\refe{0509a} and using~\refe{0509b},
we get
  $$ |A_0| \le a+1-2H < a+1-(a-2) = 3. $$
Hence, $|A_0|=2$; that is, $A_0=\{0,a\}$. Moreover, $A_1$ is contained in a
union of $|\tA_1|=2$ residue classes modulo $a$. Therefore, from $a\ge 3$,
  $$ n = |A_0|+|A_1| \le 2 + 2\lpr\frac{l-u-1}a+1\rpr
                        \le 4 + \frac23\,(l-u-1); $$
that is, $u\le l-\frac32\,n+5$, contradicting~\refe{ul1}.

This completes the proof of Theorem~\reft{stab}.

\section*{Acknowledgement}
The author is grateful to the anonymous referee for the very careful reading
of the manuscript and a number of remarks and suggestions.

\bigskip

\end{document}